\documentclass[leqno,11pt]{article}
\usepackage[spanish]{babel}
\usepackage[utf8]{inputenc}
\usepackage{amscd}
\usepackage{amsthm}
\usepackage{amssymb}
\usepackage{amsmath}
\usepackage{graphics}
\usepackage{graphicx}
\usepackage{verbatim}
\usepackage{color}

\newtheorem{mytheo}{Theorem}

\decimalpoint

\begin{document}

\begin{center}
	{\sc A Note on Consistency of the Bayes Estimator of the Density}\vspace{2ex}\\
	A.G. Nogales\vspace{2ex}\\
	Dpto. de Matem\'aticas, Universidad de Extremadura\\
	Avda. de Elvas, s/n, 06006--Badajoz, SPAIN.\\
	e-mail: nogales@unex.es
\end{center}
\vspace{.4cm}
\begin{quote}
	\hspace{\parindent} {\small {\sc Abstract.}  		
	Under mild conditions, it is shown the strong consistency of the Bayes estimator of the density. Moreover, the Bayes risk (for some common loss functions) of the Bayes estimator of the density (i.e. the posterior predictive density) reaches zero when the sample size goes to $\infty$.  In passing, a similar result is obtained for the estimation of the sampling distribution. 
		 }
\end{quote}

\vfill
\begin{itemize}
	\item[] \hspace*{-1cm} {\em AMS Subject Class.} (2010): {\em Primary\/} 62G07, 62G20
	{\em Secondary\/} 62F15
	\item[] \hspace*{-1cm} {\em Key words and phrases: } Bayesian density estimation, Bayesian estimation of the sampling distribution, posterior predictive distribution, consistency of the Bayes estimator.

\end{itemize}

\newpage

\section{Introduction}

The expression \emph{the probability of an event $A$} is an abuse of language in a statistical context as it depends on the unknown parameter (we note it $P_\theta(A)$). Before performing the experiment, this expression can be assigned a natural meaning from a Bayesian point of view as the prior predictive probability of $A$ since it is the prior mean of the probabilities   $P_\theta(A)$. But according to Bayesian philosophy, after conducting the experiment, a more appropriate estimation of  $P_\theta(A)$ is the posterior predictive probability given $\omega$ of $A$ when $\omega$ is the observed value. It has been previously proved by the author (see Nogales (2021)) that not only is this the Bayes estimator of $P_\theta(A)$ but that the posterior predictive distribution
(resp. the posterior predictive density)  is the Bayes estimator of the sampling distribution $P_\theta$ (resp. the density $p_\theta$) for the squared variation total (resp. the squared $L^1$) loss function  in the product experiment $(\Omega^n,\mathcal A^n,\{P_\theta^n\colon \theta\in(\Theta,\mathcal T,Q)\})$ corresponding to a $n$-sized sample of the unknown distribution. It should be noted that the loss functions considered derive in a natural way from the commonly used squared error loss function when estimating a real function of the parameter. 

We refer to the references within Nogales (2021) for other statistical uses of the posterior predictive distribution and some useful ways to calculate it. This paper also includes an Appendix that may be useful for  readers who are not very familiar with the mathematical notations that will be used here.

In this paper we explore the asymptotic behaviour of the posterior predictive density as Bayes estimator of the density, showing its strong consistency and that the Bayes risk goes to 0 when $n$ grows. 

\section{The framework}

Let
$$(\Omega,\mathcal A,\{P_\theta\colon \theta\in(\Theta,\mathcal T,Q)\})
$$ 
be a Bayesian experiment (where $Q$ denotes de prior distribution on the parameter space $(\Theta,\mathcal T)$) and consider the infinite product Bayesian experiment
$$(\Omega^{\mathbb N},\mathcal A^{\mathbb N},\{P_\theta^{\mathbb N}\colon \theta\in(\Theta,\mathcal T,Q)\})
$$
corresponding to an infinite sample  of the unknown distribution $P_\theta$. Let us write
$$I(\omega,\theta):=\omega,\quad J(\omega,\theta):=\theta,\quad I_n(\omega,\theta):=\omega_n\quad \mbox{and}\quad I_{(n)}(\omega):=\omega_{(n)}:=(\omega_1,\dots,\omega_n)
$$
for an integer  $n$. 

We suppose that  $P^{\mathbb N}(\theta,A):=P_\theta^{\mathbb N}(A)$ is a  Markov kernel. 
Let
$$\Pi_{\mathbb N}:=P^{\mathbb N}\otimes Q
$$
the joint distribution of the parameter and the observations, i.e., 
$$\Pi_{\mathbb N}(A\times T)=\int_TP_\theta^{\mathbb N}(A)dQ(\theta),\quad A\in\mathcal A,\; T\in\mathcal T.
$$

Being $Q:=\Pi_{\mathbb N}^{J}$ (i.e. the probability distribution of $J$ with respect to $\Pi_{\mathbb N}$),  $P_\theta^{\mathbb N}$ is a version of the conditional distribution (regular conditional probability) $\Pi_{\mathbb N}^{I|J=\theta}$. Analogously, $P_\theta^n$ is a version of the conditional distribution $\Pi_{\mathbb N}^{I_{(n)}|J=\theta}$.

Let $\beta_{Q,\mathbb N}^*:=\Pi_{\mathbb N}^{I}$, the prior predictive distribution in $\Omega^{\mathbb N}$ (so that $\beta_{Q,\mathbb N}^*(A)$ is the prior mean of the probabilities $P_\theta^{\mathbb N}(A)$). Similarly we write  $\beta_{Q,n}^*:=\Pi_{\mathbb N}^{I_{(n)}}$ for the prior predictive distribution in $\Omega^n$. So, the posterior distribution is $P^*_{\omega,\mathbb N}:=\Pi_{\mathbb N}^{J|I=\omega}$ given $\omega\in\Omega^{\mathbb N}$ satisfies 
$$\Pi_{\mathbb N}(A\times T)=\int_TP_\theta^{\mathbb N}(A)dQ(\theta)=\int_AP^*_{\omega,\mathbb N}(T)d\beta_{Q,\mathbb N}^*(\omega),\quad A\in\mathcal A,\; T\in\mathcal T.
$$
Denote
$P^*_{\omega_{(n)},n}:=\Pi_{\mathbb N}^{J|I_{(n)}=\omega_{(n)}}$ for $\omega_{(n)}\in\Omega^n$, the posterior distribution given $\omega_{(n)}\in\Omega^n$.

Write ${P^*_{\omega_{(n)},n}}^{\!\!\!\!\!\!\!\!\!\!P\,\,\,\,\,}$ for the posterior predictive distribution given $\omega_{(n)}\in\Omega^n$ defined for $A\in\mathcal A$ as
$${P^*_{\omega_{(n)},n}}^{\!\!\!\!\!\!\!\!\!\!P\,\,\,\,\,}(A)=
\int_\Theta P_\theta(A)dP^*_{\omega_{(n)},n}(\theta).
$$
So ${P^*_{\omega_{(n)},n}}^{\!\!\!\!\!\!\!\!\!\!P\,\,\,\,\,}(A)$ is nothing but the posterior mean given $\omega_{(n)}\in\Omega^n$ of the probabilities $P_\theta(A)$.

In the dominated case we can assume without loss of generality that the dominating measure $\mu$ is a probability measure (because of (1) below). We write $p_\theta=dP_\theta/d\mu$. The likelihood function $\mathcal L(\omega,\theta):=p_\theta(\omega)$ is supposed $\mathcal A\times\mathcal T$-measurable. 

We have that, far all $n$ and every event $A\in\mathcal A$, 
\begin{gather*}\begin{split}
		{P^*_{\omega_{(n)},n}}^{\!\!\!\!\!\!\!\!\!\!P\,\,\,\,\,}(A)&=
		\int_\Theta P_\theta(A)dP^*_{\omega_{(n)},n}(\theta)=
		\int_\Theta \int_A p_\theta(\omega')d\mu(\omega')dP^*_{\omega_{(n)},n}(\theta)
			\\&	=
		\int_A\int_\Theta  p_\theta(\omega')dP^*_{\omega_{(n)},n}(\theta)d\mu(\omega')
	\end{split}\end{gather*}
which proves that
$${p^*_{\omega_{(n)},n}}^{\!\!\!\!\!\!\!\!\!\!P\,\,\,\,\,}(\omega'):=\int_{\Theta} p_\theta(\omega')dP^*_{\omega_{(n)},n}(\theta)
$$
is a $\mu$-density of  ${P^*_{\omega_{(n)},n}}^{\!\!\!\!\!\!\!\!\!\!P\,\,\,\,\,}$ that we recognize as the posterior predictive density on $\Omega$ given $\omega_{(n)}$.

In the same way
$${p^*_{\omega,\mathbb N}}^{\!\!\!\!P}(\omega'):=\int_{\Theta} p_\theta(\omega')dP^*_{\omega,\mathbb N}(\theta)
$$
is a $\mu$-density of ${P^*_{\omega,\mathbb N}}^{\!\!\!\!P}$, the posterior predictive density on $\Omega$ given $\omega\in\Omega^{\mathbb N}$. 

In the next we will assume the following additional regularity conditions: 
\begin{itemize}
	\item[(i)] $(\Omega,\mathcal A)$ is a Borel standard space, 
	
	\item[(ii)] $\Theta$ is a Borel subset of a Polish space and  $\mathcal T$ is its Borel $\sigma$-field, and
	\item[(iii)] $\{P_\theta\colon \theta\in\Theta\}$ is identifiable.
\end{itemize}

According to Nogales (2021), the posterior predictive distribution
${P^*_{\omega_{(n)},n}}^{\!\!\!\!\!\!\!\!\!\!P\,\,\,\,\,}$
(resp. the posterior predictive density ${p^*_{\omega_{(n)},n}}^{\!\!\!\!\!\!\!\!\!\!P\,\,\,\,\,}$) is the Bayes estimator of the sampling distribution $P_\theta$ (resp. the density $p_\theta$) for the squared variation total (resp. the squared $L^1$) loss function  in the product experiment $(\Omega^n,\mathcal A^n,\{P_\theta^n\colon \theta\in(\Theta,\mathcal T,Q)\})$. Analogously, the posterior predictive distribution 
${P^*_{\omega,\mathbb N}}^{\!\!\!\!P}$ (resp. the posterior predictive density ${p^*_{\omega,\mathbb N}}^{\!\!\!\!P}$) is the Bayes estimator of the sampling distribution $P_\theta$ (resp. the density $p_\theta$) for the squared variation total (resp. the squared $L^1$) loss function  in the product experiment $(\Omega^{\mathbb N},\mathcal A^{\mathbb N},\{P_\theta^{\mathbb N}\colon \theta\in(\Theta,\mathcal T,Q)\})$.

As a particular case of a well known result about the total variation distance between two probability measures and the $L^1$-distance between their densities, we have that
$$\sup_{A\in\mathcal A}\left|{P^*_{\omega_{(n)},n}}^{\!\!\!\!\!\!\!\!\!\!P\,\,\,\,\,}(A)-P_\theta(A)\right|=\frac12\int_{\Omega}\left|{p^*_{\omega_{(n)},n}}^{\!\!\!\!\!\!\!\!\!\!P\,\,\,\,\,}-p_\theta\right|d\mu.\qquad (1)
$$

\section{The main result}

We wonder if the Bayes risk of the Bayes estimator ${P^*_{\omega_{(n)},n}}^{\!\!\!\!\!\!\!\!\!\!P\,\,\,\,\,}$ of the sampling distribution $P_\theta$ goes to zero when  $n\to\infty$, i.e., if 
$$\lim_n\int_{\Omega^{\mathbb N}\times\Theta}\sup_{A\in\mathcal A}\left|{P^*_{\omega_{(n)},n}}^{\!\!\!\!\!\!\!\!P\,\,\,\,}(A)-P_\theta(A)\right| d\Pi_{\mathbb N}(\omega,\theta)=0.
$$

In terms of densities, the question is whether  the Bayes risk of the Bayes estimator  ${p^*_{\omega_{(n)},n}}^{\!\!\!\!\!\!\!\!\!\!P\,\,\,\,\,}$ of the density $p_\theta$ goes to zero when  $n\to\infty$, i.e., if 
$$\lim_n \int_{\Omega^{\mathbb N}\times\Theta}\left(\int_{\Omega}\left|{p^*_{\omega_{(n)},n}}^{\!\!\!\!\!\!\!\!\!\!P\,\,\,\,\,}(\omega')-p_\theta(\omega')\right|d\mu(\omega')\right)^2d\Pi_{\mathbb N}(\omega,\theta)=0.
$$

Let us consider the auxiliary Bayesian experiment
$$(\Omega\times\Omega^{\mathbb N},\mathcal A\times\mathcal A^{\mathbb N},\{\mu\times P_\theta^{\mathbb N}\colon \theta\in(\Theta,\mathcal T,Q)\}).
$$

For $\omega'\in\Omega$, $\omega\in\Omega^n$ and $\theta\in\Theta$, we will continue to denote $I(\omega',\omega,\theta)=\omega$ and $J(\omega',\omega,\theta)=\theta$, and now we write  $I'(\omega',\omega,\theta)=\omega'$.

The new prior predictive distribution is $\mu\times \beta_{Q,n}^*$ since
$$(\mu\times\Pi_{\mathbb N})^{(I',I_{(n)})}(A'\times A_{(n)})=\mu(A')\cdot \beta_{Q,n}^*(A_{(n)})=(\mu\times \beta_{Q,n}^*)(A'\times A_{(n)}).
$$
To compute the new posterior distributions, notice that 
\begin{gather*}
	(\mu\times\Pi_{\mathbb N})(A'\times I_{(n)}^{-1}( A_{(n)} ) \times T)=\\
	\int_{A'\times I_{(n)}^{-1}( A_{(n)} ) } 
	(\mu\times\Pi_{\mathbb N})^{J| (I',I_{(n)}) =
		(\omega',\omega_{(n)}) } (T) d(\mu\times\Pi_{\mathbb N})^{(I',I_{(n)})} (\omega',\omega_{(n)}).
\end{gather*}
On the other hand,
\begin{gather*}
	(\mu\times\Pi_{\mathbb N})(A'\times I_{(n)}^{-1}( A_{(n)} ) \times T)=\mu(A')\cdot
	\Pi_{\mathbb N}(I_{(n)}^{-1}( A_{(n)} ) \times T)=\\
	\mu(A')\cdot\int_{A_{(n)}}
	P^*_{\omega_{(n)},n}(T)d\beta^*_{Q,n}(\omega_{(n)})
	=\int_{A'\times A_{(n)}}
	P^*_{\omega_{(n)},n}(T)d(\mu\times\beta^*_{Q,n})(\omega',\omega_{(n)}).
\end{gather*}
So,
$$P^*_{\omega_{(n)},n}=(\mu\times\Pi_{\mathbb N})^{J|(I',I_{(n)})=(\omega',\omega_{(n)})}.
$$
It follows that if $f\in L^1(Q)$ then 
$$E_{P^*_{\omega_{(n)},n}}(f)=E_{\mu\times\Pi_{\mathbb N}}[f\mid(I',I_{(n)})=(\omega',\omega_{(n)})].
$$

When $\mathcal A'_{(n)}:=(I',I_{(n)})^{-1}(\mathcal A\times\mathcal A^n)$, 
we have that  $(\mathcal A'_{(n)})_n$ is an increasing sequence of  sub-$\sigma$-fields of $\mathcal A\times\mathcal A^{\mathbb N}$ such that  $\mathcal A\times\mathcal A^{\mathbb N}=\sigma(\cup_n\mathcal A'_{(n)})$. According to the martingale convergence theorem of Lévy, if $Y$ es $(\mathcal A\times\mathcal A^{\mathbb N}\times\mathcal T)$-measurable and $\mu\times\Pi_{\mathbb N}$-integrable, then 
$$E_{\mu\times \Pi_{\mathbb N}}(Y|\mathcal A'_{(n)})
$$
converges $(\mu\times\Pi_{\mathbb N})$-a.e. and in $L^1(\mu\times\Pi_{\mathbb N})$ to $Y=E_{\mu\times\Pi_{\mathbb N}}(Y|\mathcal A'\times\mathcal A^{\mathbb N})$. 

Let us consider the $\mu\times\Pi_{\mathbb N}$-integrable function 
$$Y(\omega',\omega,\theta):=p_\theta(\omega'). 
$$
Next we show that
$${p^*_{\omega,\mathbb N}}^{\!\!\!\!P}(\omega')=E_{\mu\times\Pi_{\mathbb N}}(Y\mid (I',I)=(\omega',\omega)).\qquad (2)
$$
Indeed, given $A',A\in\mathcal A$, we have that
\begin{gather*}\begin{split}
		\int_{(I',I)^{-1}(A'\times A)} & p_\theta(\omega') d(\mu\times\Pi_{\mathbb N})(\omega',\omega,\theta)=\int_A\int_\Theta\int_{A'}p_\theta(\omega')d\mu(\omega')dP^*_{\omega,\mathbb N}(\theta)d\beta^*_{Q,\mathbb N}(\omega)\\
		&=\int_A\int_\Theta P_\theta(A')dP^*_{\omega,\mathbb N}(\theta)d\beta^*_{Q,\mathbb N}(\omega)
		=\int_A {P^*_{\omega,\mathbb N}}^{\!\!\!\!P}(A')d\beta^*_{Q,\mathbb N}(\omega)\\
		&=\int_{A'}\int_A {p^*_{\omega,\mathbb N}}^{\!\!\!\!P}(\omega')d\mu(\omega')d\beta^*_{Q,\mathbb N}(\omega)=
		\int_{A'\times A} {p^*_{\omega,\mathbb N}}^{\!\!\!\!P}(\omega')
		d(\mu\times\Pi_{\mathbb N})^{(I',I)}(\omega',\omega)
\end{split}\end{gather*}
which proves (2). 

Analogously it can be shown that 
$${p^*_{\omega_{(n)},n}}^{\!\!\!\!\!\!\!\!P\,\,\,\,}(\omega')=E_{\mu\times\Pi_{\mathbb N}}(Y\mid (I',I_{(n)})=(\omega',\omega_{(n)})).\qquad (3)
$$

Hence, it follows from the aforementioned theorem of Lévy that
$$\lim_n {p^*_{\omega_{(n)},n}}^{\!\!\!\!\!\!\!\!P\,\,\,\,}(\omega') = {p^*_{\omega,\mathbb N}}^{\!\!\!\!P}(\omega'),\quad (\mu\times \Pi_{\mathbb N})-\hbox{a.e.} \qquad (4)
$$
and
$$\lim_n\int_{\Omega\times\Omega^{\mathbb N}\times\Theta}\left|{p^*_{\omega_{(n)},n}}^{\!\!\!\!\!\!\!\!P\,\,\,\,}(\omega')-{p^*_{\omega,\mathbb N}}^{\!\!\!\!P}(\omega')\right| d(\mu\times\Pi_{\mathbb N})(\omega',\omega,\theta)=0,
$$
i.e., 
$$\lim_n\int_{\Omega^{\mathbb N}\times\Theta}\int_\Omega\left|{p^*_{\omega_{(n)},n}}^{\!\!\!\!\!\!\!\!P\,\,\,\,}(\omega')-{p^*_{\omega,\mathbb N}}^{\!\!\!\!P}(\omega')\right| d\mu(\omega')d\Pi_{\mathbb N}(\omega,\theta)=0.\qquad (5)
$$

On the other hand, as a consequence of a known theorem of Doob (see Theorem 6.9 and Proposition 6.10 from Ghosal et al. (2017, p. 129, 130)) we have that, for every  $\omega'\in\Omega$,
$$\lim_n \int_\Theta p_{\theta'}(\omega')dP^*_{\omega_{(n)},n}(\theta')=p_\theta(\omega'),\quad P_\theta^{\mathbb N}-\hbox{a.e.}
$$
for $Q$-almost every $\theta$.

Hence
$$\lim_n {p^*_{\omega_{(n)},n}}^{\!\!\!\!\!\!\!\!\!\!P\,\,\,\,\,}(\omega')=p_\theta(\omega'),\quad P_\theta^{\mathbb N}-\hbox{a.e.}
$$
for $Q$-almost every $\theta$, i.e., given $\omega'\in\Omega$ there exists $T_{\omega'}\in\mathcal T$ such that $Q(T_{\omega'})=0$ and, $\forall \theta\notin T_{\omega'}$,
$$\lim_n {p^*_{\omega_{(n)},n}}^{\!\!\!\!\!\!\!\!\!\!P\,\,\,\,\,}(\omega')=p_\theta(\omega'),\quad P_\theta^{\mathbb N}-\hbox{a.e.}.
$$
So, for $\theta\notin T_{\omega'}$, there exists $N_{\theta,\omega'}\in\mathcal A^{\mathbb N}$ such that  $P_\theta^{\mathbb N}(N_{\theta,\omega'})=0$ and
$$\lim_n {p^*_{\omega_{(n)},n}}^{\!\!\!\!\!\!\!\!\!\!P\,\,\,\,\,}(\omega')=p_\theta(\omega'),\quad\forall\omega\notin N_{\theta,\omega'},\forall \theta\notin T_{\omega'}, \forall \omega'\in\Omega.
$$
In particular, 
$$\lim_n {p^*_{\omega_{(n)},n}}^{\!\!\!\!\!\!\!\!\!\!P\,\,\,\,\,}(\omega')=p_\theta(\omega'),\quad \mu\times P_\theta^{\mathbb N}-\hbox{a.e.}\qquad (6)
$$

From (4) and (6) it follows that  $p_\theta(\omega')={p^*_{\omega,\mathbb N}}^{\!\!\!\!\!\!P\,\,\,\,\,}(\omega')$, $\mu\times P_\theta^{\mathbb N}-\hbox{a.e.}$

From this and (5) it follows that 
$$\lim_n\int_{\Omega^{\mathbb N}\times\Theta}\int_\Omega\left|{p^*_{\omega_{(n)},n}}^{\!\!\!\!\!\!\!\!P\,\,\,\,}(\omega')-p_\theta(\omega')\right| d\mu(\omega')d\Pi_{\mathbb N}(\omega,\theta)=0, 
$$
i.e., the risk of the Bayes estimator of the density for the $L^1$ loss function goes to 0 when $n\to\infty$.  

It follows from this and (1) that 
$$\lim_n\int_{\Omega^{\mathbb N}\times\Theta}\sup_{A\in\mathcal A}\left|{P^*_{\omega_{(n)},n}}^{\!\!\!\!\!\!\!\!P\,\,\,\,}(A)-P_\theta(A)\right| d\Pi_{\mathbb N}(\omega,\theta)=0, 
$$
i.e., the risk of the Bayes estimator of the sampling distribution $P_\theta$ for the variation total loss function goes to 0 when $n\to\infty$.

We wonder if these results remain true for the squared versions of the loss functions. The answer is yes because of the following general result: let $(X_n)$ be a sequence of r.r.v. on a probability space $(\Omega,\mathcal A,P)$ such that $\lim_n\int|X_n|dP=0$. If there exists $a>0$ such that $|X_n|\le a$, for all $n$, then $\lim_n\int|X_n|^2dP=0$ because 
$$0\le\int|X_n|^2dP\le a\int|X_n|dP\rightarrow_n 0.
$$
In our case $a=2$, $P:=\Pi_{\mathbb N}$ and
$$X_n:=\int_\Omega\left|{p^*_{\omega_{(n)},n}}^{\!\!\!\!\!\!\!\!P\,\,\,\,}(\omega')-{p^*_{\omega,\mathbb N}}^{\!\!\!\!P}(\omega')\right| d\mu(\omega'), \quad\hbox{or}\quad 
X_n:=\sup_{A\in\mathcal A}\left|{P^*_{\omega_{(n)},n}}^{\!\!\!\!\!\!\!\!P\,\,\,\,}(A)-P_\theta(A)\right|.
$$

So, we have proved the following result. 

\begin{mytheo}\rm Let $(\Omega,\mathcal A,\{P_\theta\colon \theta\in(\Theta,\mathcal T,Q)\})$ be a Bayesian experiment dominated by a $\sigma$-finite measure $\mu$. Let us suppose that $(\Omega,\mathcal A)$ is a  Borel standar space, that $\Theta$ is a Borel subset of a Polish space and $\mathcal T$ is its Borel  $\sigma$-field. Suppose also that the likelihood function  $\mathcal L(\omega,\theta):=p_\theta(\omega)=\frac{dP_\theta}{d\mu}(\omega)$ is $\mathcal A\times\mathcal T$-measurable and the family $\{P_\theta\colon \theta\in\Theta\}$ is identifiable. Then:
	\begin{itemize}
		\item[(a)] The posterior predictive density  ${p^*_{\omega_{(n)},n}}^{\!\!\!\!\!\!\!\!P\,\,\,\,}$ 
		is the Bayes estimator of the density  $p_\theta$ in the product experiment  $(\Omega^n,\mathcal A^n,\{P_\theta^n\colon \theta\in(\Theta,\mathcal T,Q)\})$ for the  squared $L^1$ loss function. Moreover the risk function converges to 0 both for the $L^1$ loss function and the squared  $L^1$ loss function. 
		
		\item[(b)] The posterior predictive distribution  ${P^*_{\omega_{(n)},n}}^{\!\!\!\!\!\!\!\!P\,\,\,\,}$ 
		is the Bayes estimator of the sampling distribution  $P_\theta$ in the product experiment  $(\Omega^n,\mathcal A^n,\{P_\theta^n\colon \theta\in(\Theta,\mathcal T,Q)\})$ for the squared variation total loss function. Moreover the risk function converges to 0 both for the variation total loss function and the squared variation total loss function. 
		
		\item[(c)] The posterior predictive density is a strongly consistent estimator of the density $p_\theta$, i.e.,
		$$\lim_n {p^*_{\omega_{(n)},n}}^{\!\!\!\!\!\!\!\!\!\!P\,\,\,\,\,}(\omega')=p_\theta(\omega'),\quad \mu\times P_\theta^{\mathbb N}-\hbox{a.e.}
		$$
		for $Q$-almost every $\theta\in\Theta$. 
		
	\end{itemize}
\end{mytheo}

\section{Acknowledgements.}
This paper has been supported by the Junta de Extremadura (Spain) under the grant Gr18016.
\vspace{1ex}

\section* {References:}

\begin{itemize}

	\item Ghosal, S., Vaart, A.v.d. (2017) Fundamentals of Nonparametric Bayesian Inference, Cambridge University Press, Cambridge UK. 
	
    \item Nogales, A.G. (2021), On Bayesian Estimation of Densities and Sampling Distributions: the Posterior Predictive Distribution as the Bayes Estimator, to appear in Statistica Neerlandica. 
	
\end{itemize}

\end{document}